*Article*

# Experimental and Theoretical Investigations of a Ground Source Heat Pump System for Water and Space Heating Applications in Kazakhstan

Yelnar Yerdesh [1,2,*], Tangnur Amanzholov [1,2], Abdurashid Aliuly [1,2], Abzal Seitov [1,2], Amankeldy Toleukhanov [2], Mohanraj Murugesan [3], Olivier Botella [4], Michel Feidt [4], Hua Sheng Wang [5], Alexandr Tsoy [6] and Yerzhan Belyayev [1,2,*]

[1] Department of Mechanics, Al-Farabi Kazakh National University, Almaty 050040, Kazakhstan
[2] Department of Mechanical Engineering, Satbayev University, Almaty 050013, Kazakhstan
[3] Department of Mechanical Engineering, Hindusthan College of Engineering and Technology, Coimbatore 641032, India
[4] Université de Lorraine, CNRS, LEMTA, F-54000 Nancy, France
[5] School of Engineering and Materials Sciences, Queen Mary University of London, London E1 4NS, UK
[6] Department of Machines and Devices of Manufacturing Processes, Almaty Technological University, Almaty 050012, Kazakhstan
* Correspondence: yelnaryerdesh@gmail.com (Y.Y.); y.belyayev@satbayev.university (Y.B.); Tel.: +7-707-803-0104 (Y.Y.); +7-771-491-3344 (Y.B.)

**Abstract:** The ground source heat pump heating system is considered as one of the best solutions for the transition towards green heating under the continental climate conditions like Kazakhstan. In this paper, experimental and theoretical investigations were carried out to develop a ground source heat pump-based heating system under the weather conditions in Kazakhstan and to evaluate its thermodynamic performance. The water-to-water heat pump heating system, integrated with a ground source heat exchanger and used refrigerant R134a, was designed to provide hot water to meet the requirements for space heating. The predicted values of the coefficient of performance and the experimental results were found to be in good agreement within 6.2%. The thermodynamic performance of the system was also assessed using various environment-friendly refrigerants, such as R152a, R450A, R513A, R1234yf and R1234ze, as potential replacements for R134a. Although R152a is found to be a good alternative for R134a in terms of coefficient of performance and total equivalent warming impact, its flammability hinders its application. The heating system using refrigerants R450A, R513A, R1234yf and R1234ze shows 2–3% lower coefficient of performance than that of R134a. The highest exergy destruction is found to be attributed to the compressor, followed by the expansion valve, evaporator, and condenser. Considering their low flammability and low environmental impact, R450A, R513A, R1234yf and R1234ze are identified as valuable replacements for R134a.

**Keywords:** ground source heat pump; experimental setup; coefficient of performance; energy efficiency; exergy efficiency; environmental-friendly refrigerants

## 1. Introduction

Heat pump systems are identified as energy efficient heating equipment used for water heating and space heating in continental climatic regions. The performance of heat pumps is improved by integrating with renewable energy sources like solar, geothermal, and solar-geothermal hybrid modes. Due to poor solar availability during winter in continental climates, the concept of integrating solar energy with heat pump systems is not feasible. Presently, the ground source heat pumps (GSHPs) are considered as the most efficient space heating system in continental climatic regions in Kazakhstan [1–4]. Low-temperature underground thermal energy source at 30–40 °C is available at a depth of 200 to 400 m [5,6]. To maintain thermal comfort of the room during winter, thermal energy is extracted from the ground; whereas, during summer, the excess thermal energy from the room is pumped and stored into the ground. This is especially valuable for Kazakhstan, whose majority of territory belongs to the continental climate, which is characterized by the presence of contrasted seasons from cold winters to hot summers [7]. Ambient air temperature significantly affects not only the building heat demand but also the shallow geothermal energy characteristics. For example, in severe winter conditions, since the building heat demand is high, the operation of a GSHP causes a quick depletion of the soil thermal energy around borehole heat exchanger (BHE) [8,9]. This leads to a gradual decrease in soil temperature. Such problems may arise for the central, eastern, and northern regions of Kazakhstan when severe winter occurs. In order to restore the thermal energy that has been depleted from the soil during the summer heating season, the excess solar thermal and waste heat from space cooling is pumped into the soil through ground heat exchangers [8–10]. Depending on the solar irradiance and the room cooling demand, this thermal energy may not be sufficient to restore the thermal energy depleted during the heating season. For the southern and western (Caspian) regions of Kazakhstan, where milder winter conditions are observed, ground thermal depletion may not be significant. In this case, when surplus solar heat and waste heat from the heat pump cooling mode in summer are pumped into the ground, not only the depleted thermal energy is restored but also excess thermal energy may be stored. Seasonal underground thermal energy storage is implemented using the borehole thermal energy storage (BTES) and aquifer thermal energy storage (ATES) methods [6,11]. In the case of BTES, underground closed-loop type BHEs are used, usually U-shaped vertical heat exchangers. In contrast to BTES system, ATES is an open-loop type system, where groundwater from shallow bore holes is circulated between two or more wells. During the heating season, the groundwater is circulated from heat source well(s) to heat sink well(s) using a pump; in summer, it is circulated in the opposite direction.

The two types of ground heat exchangers integrated with heat pump system are: (a) horizontal, and (b) vertical types. The horizontal ground heat exchangers (GHEs) are buried in trenches with a depth of 1–2 m, where the distance between the trenches is 2–3 m and the trench length is 24–48 m [12]. Due to low cost and easy installation, such systems are popular in warm regions for building heating and cooling. The horizontal type GHEs are recommended in the southern and western parts of Kazakhstan. In these locations, the higher atmospheric temperature and more solar irradiance availability positively influence the utilization of horizontal GHEs. A horizontal GHE is a closed-loop type system, where a heat transfer fluid (HTF) circulates inside to exchange heat only with the ambient soil. However, the horizontal heat exchanger is not suitable for the extreme winter conditions with low solar irradiance [13,14].

The most common type of GHE is a vertical type that can be used on a limited surface area without disturbing the landscape. The BHEs and geothermal energy piles (GEPs) [1] can harvest natural ground temperature, which is used as a heat source for heating in winter and a heat sink for cooling in summer. The GEPs serve as load-bearing deep foundations for buildings as well as GHEs for GSHPs. The most commonly used BHEs are U-shaped pipes made of high-density polyethylene (HDPE), which are installed in a well (of various diameters) backfilled with cement/bentonite grout. These vertical boreholes are buried at a depth of about 200 m, which corresponds to a shallow geothermal energy. Depending on the soil or rock type, the initial investment for drilling wells is always high. To reduce costs while maintaining the thermal potential of the soil, it is recommended to drill several wells down to 100 m. A spacing of 5–8 m is required to minimize thermal interference between several BHEs [15,16]. The BHE (GEP) is a closed-loop type GHE, where HTF circulates heat between heat pump evaporator and the ground. The BHEs are integrated with a water-to-water heat pump for efficient use of the shallow geothermal energy. Such U-shaped BHE-based GSHPs are recommended to be used in almost the entire territory of Kazakhstan, especially where severe winters are observed. The influence of the atmospheric air temperature on vertical GHEs is found to be insignificant [17–20].

Another type of GHEs with a closed-loop structure are shallow ground basket-type heat exchangers, which are separate pipe loops arranged in a spiral or in a special twisted shape, buried at a shallow depth of 1–6 m [21,22]. These type of shallow subsurface GHEs, also called energy baskets, can be immersed in a water body such as a dam or a lake away from the building [21–25]. Like horizontal GHE based GSHPs, these energy baskets significantly reduce the capital cost and are suitable for mild winter conditions or for buildings near lakes/dams. Such geothermal heat exchangers are not popular for the continental climatic conditions of Kazakhstan.

Closed-loop systems have been discussed above. In contrast to these systems, where only heat transfer occurs, both heat and mass transfer occurs in GHEs with an open-loop structure. An open-water loop system involves underground water that is being abstracted from shallow geothermal aquifer layers, and then returned after heat is extracted or added. Therefore, such a system consists of at least two wells where underground water from the aquifer is taken and reinjected from one well to another. One well is used for heat storage and the other for cold storage, according to the season. The wells are located at several tens of meters from each other, excluding the mutual influence of warm and cold "bells", and are connected on the ground to a water-to-water type heat pump by a pipeline [26]. In most systems, the winter temperature of the reinjected water is 6–9 °C, while the summer temperature is 15–25 °C [26]. One of the important advantages of such a system is the large volumes of water abstracted, which results in a large heat exchange capacity and can be used in large scale applications [27–30].

Among these types of GSHPs, heat pumps with BHEs and wells from aquifers are the most mature technologies for utilizing shallow geothermal energy. The BHEs do not require special preliminary survey of geological or hydrogeological structure on the presence of aquifers and are applicable to all climates. The seasonal variation of the ground temperature due to seasonal change in the ambient air temperature only occurs at the depth of 15 m below ground surface for cold regions such as Kazakhstan, Russia, Northern China, Mongolia, Northern USA and Canada. Therefore, the use of such vertical BHE (below 50–150 m depth) as part of the GSHP will be efficient and widespread. In this paper, an experimental and theoretical study of a GSHP is carried out using U-shape BHEs under the weather conditions of the city of Almaty, in the Republic of Kazakhstan.

The vapor compression water-to-water heat pump is used to enhance the low-grade heat extracted through GHEs from shallow geothermal resources. The vapor-compression heat pump cycle can be described by the reverse Carnot cycle, an ideal circular process consisting of two adiabatic and two isothermal processes [31,32]. The main working fluid is a refrigerant, where the heat exchange process takes place with phase change. In an indirect expansion system, a water-to-water or brine-to-water heat pump is a main element [33,34]. GSHPs with aquifers as heat sources use a water-to-water heat pump. BHE-based GSHPs use both water-to-water and brine-to-water heat pumps, depending on the source and sink temperature ranges. Recent works in this area are reported in [35–37].

Many research initiatives have been devoted to heat pump systems with vertical and horizontal ground heat exchangers under various weather climates [2,3,7,8,19]. However, to the best of our knowledge, thermodynamic performances for the specific weather conditions of Kazakhstan have not yet been reported. Moreover, halogen-based refrigerants like R134a, R407C and R410A have been widely used in heat pump systems integrating with ground heat exchangers. These refrigerants have a high GWP, and many researchers are looking for low GWP alternatives [35–37]. Hence, the possibility of using ground source heat pumps using environment friendly refrigerant mixtures for hot water and space heating applications in Kazakhstan is considered in this research. A detailed energy and exergy analysis is carried out to assess both energy conversion and energy losses in the system.

## 2. Apparatus

The experimental setup installed in Almaty is depicted in Figure 1. Two vertical U-shaped BHEs are used to extract the low-grade ground thermal energy using a single U-pipe and double U-pipe of 50 m depth. A water-to-water heat pump with a rated heating capacity of 5 kW is developed and integrated with BHEs for water heating and space heating. A hot water storage tank of 300 L is used to store the thermal energy.

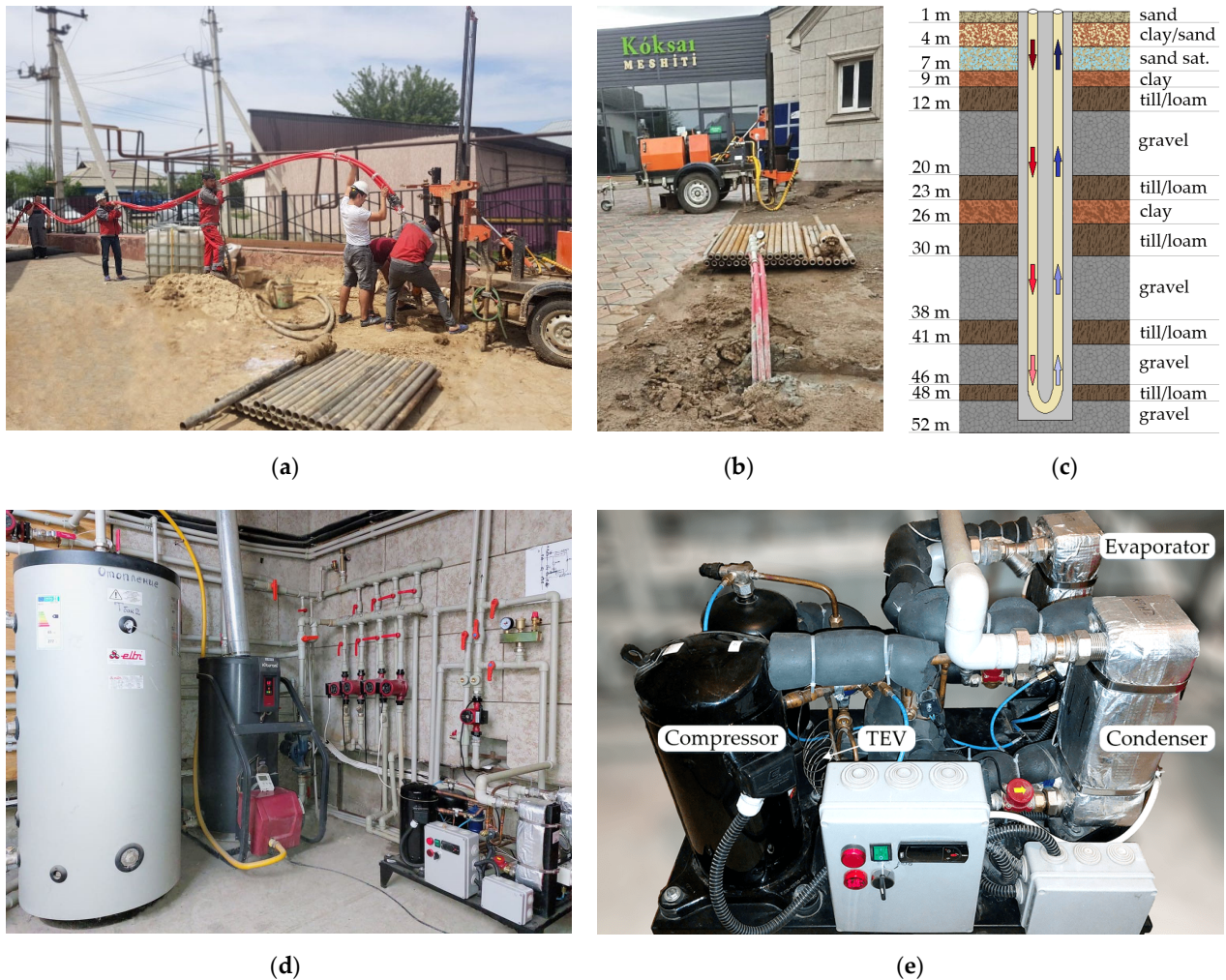

**Figure 1.** Photographs of the experimental facility in Koksai Mosque, in the Almaty region of Kazakhstan: (**a**) BHE installation; (**b**) Buried BHE; (**c**) Ground structure; (**d**) Heat delivery unit; (**e**) Water-to-water heat pump.

To install the two BHEs, two wells of depth 50 m have been drilled. High-density polyethylene (HDPE) piping material of 32 mm diameter is used for the U-shaped heat exchangers. After immersing the heat exchangers into the well, the voids are filled with grout, which is a solution of water (70%), cement (24%) and bentonite (6%). After installing the GHEs, pressure drops and air in the hydraulic loop are checked. A thermal response test (TRT) is carried out to determine the soil thermal conductivity and the BHE thermal resistance. A water-to-water heat pump operates for two modes (heating and cooling) using a reversing valve. The performance of the water-to-water heat pump, integrated with ground heat exchanger operating in heating mode, is tested for water heating and space heating.

Figure 2 shows the schematic of the water-to-water heat pump integrated with GHE. The system has three loops: GHE loop, heat pump refrigerant loop and hot water storage loop. The GHE circuit consists of a U-shaped BHE, a circulation pump, an expansion tank, and an evaporator of the heat pump. The heat pump refrigerant loop consists of a compressor, a condenser, an expansion device, and an evaporator. The detailed specifications of the components of the water-to-water heat pump are given in Tables 1–3. Accessories such as filter-drier, liquid receiver and sight glass are installed in the liquid line. The reversing valve is used to divert the refrigerant flow according to the requirement. The hot water storage tank loop consists of a condenser, a circulation pump, a 300-L storage tank and an expansion tank. The temperature and pressure of refrigerants in the water-to-water heat pump loop are measured with thermocouples and pressure gauges (see Figure 2).

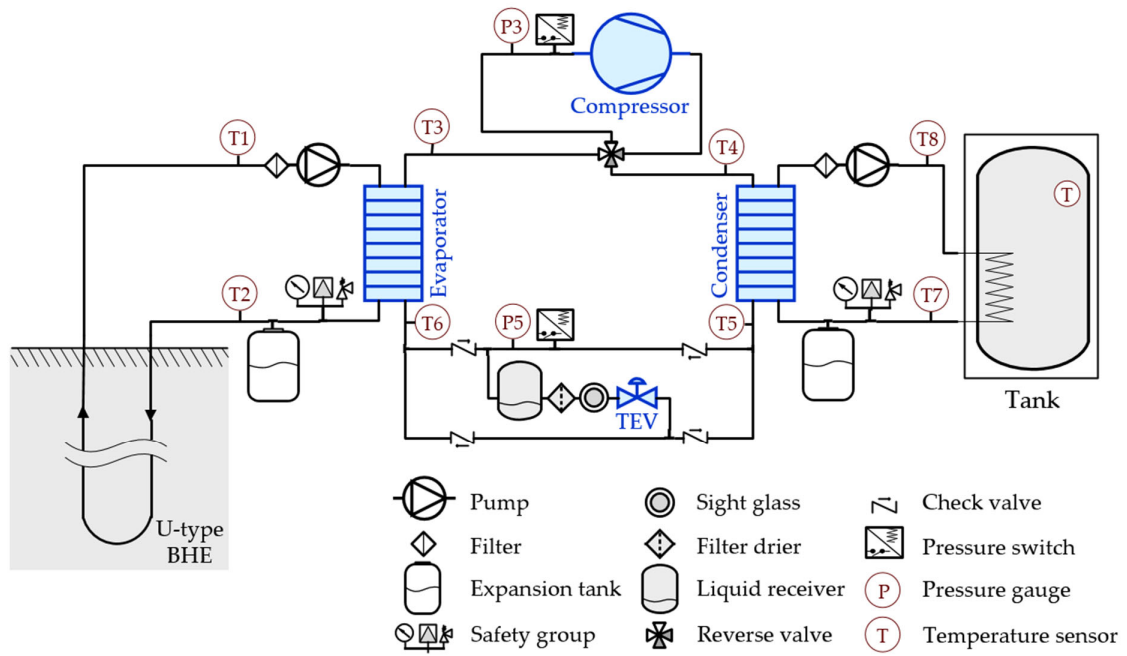

**Figure 2.** Schematic of the experimental setup.

**Table 1.** Compressor parameters supplied by the manufacturer.

| Parameter | Unit | Compressor |
|---|---|---|
| Model | - | Copeland ZR28K3E-PFJ |
| Displacement @ 50 Hz | m³/h | 6.83 |
| Max. Internal Free Volume | L | 2.90 |
| Dimensions (L × W × H) | mm | 246 × 246 × 377 |
| Oil type | - | POE RL32-3MAF |
| Oil quantity | L | 1.12 |
| Heating capacity * | kW | 4.83 |
| Input power * | kW | 1.38 |
| *COP* * | - | 2.57 |
| Mass flow rate * | g/s | 24.7 |

* Parameter obtained at $T_{evap}$ = 0 °C and $T_{cond}$ = 50 °C, $T_{SH}$ = 10 °C, $T_{SC}$ = 5 °C.

**Table 2.** Expansion device parameters supplied by manufacturer.

| Parameter | Unit | TEV |
|---|---|---|
| Model | - | Danfoss T2–4 |
| Maximum evaporating temperature | °C | 15 |
| Minimum evaporating temperature | °C | −45 |
| Max. working pressure | kPa | 3400 |
| Pressure drop * | kPa | 1024 |
| Nominal capacity * | kW | 6.4 |
| Min. capacity * | kW | 1.6 |

* Parameters obtained at $T_{evap} = 0$ °C and $T_{cond} = 50$ °C, $T_{SH} = 10$ °C, $T_{SC} = 5$ °C.

**Table 3.** Heat exchanger parameters supplied by manufacturer.

| Parameter | Unit | Evaporator | | Condenser | |
|---|---|---|---|---|---|
| Model | - | B3-027-20-H | | B3-027-30-H | |
| Heat transfer area | m² | 0.468 | | 0.728 | |
| Dimensions (L × W × H) | mm | 111 × 54 × 311 | | 111 × 81 × 311 | |
| Weight | kg | 3.8 | | 5.1 | |
| Material | - | AISI 316L | | AISI 316L | |
| Flow arrangement | - | Counterflow | | Counterflow | |
| Fluid side | - | R134a | Water | R134a | Water |
| Number of plates | - | 10 | 9 | 15 | 14 |
| Heat transfer coefficient * | W/(m²·K) | 2177.1 | 5224.0 | 1347.1 | 7778.2 |
| Total pressure drop * | kPa | 10.47 | 1.89 | 0.21 | 2.35 |

* Parameters obtained at $T_{evap} = 0$ °C and $T_{cond} = 50$ °C, $T_{SH} = 10$ °C, $T_{SC} = 5$ °C.

The present work reports the experimental data of a GSHP heating system and the thermodynamic analysis based on the 1st and 2nd laws of thermodynamics [31,32]. Figure 3 shows the vapor-compression heat pump cycle on the pressure-enthalpy diagram. The HTF leaving the BHE enters the evaporator at point 1, and the heat is transferred from the HTF to the refrigerant. The heat transfer of the HTF in both the BHE and the evaporator is single-phase convective heat transfer. The HTF is heated in the condenser from point 7 to 8 and then transfers heat to water in the thermal energy storage tank. The heat transfer of the HTF in both the condenser and the tank is also single-phase convective heat transfer.

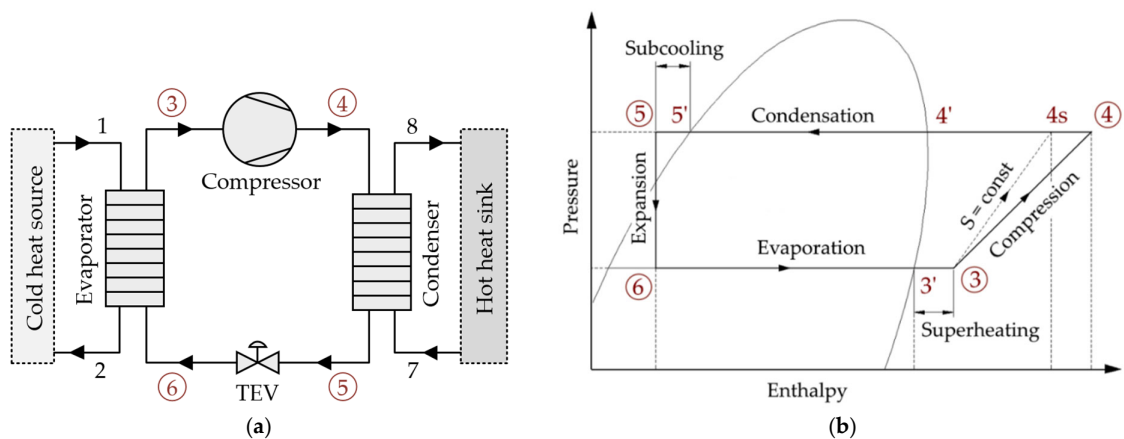

**Figure 3.** Schematic representation of the thermodynamic model: (**a**) Water-to-water heat pump; (**b**) Pressure-enthalpy diagram.

The processes 3–4, 4–5, 5–6, and 6–3 represent the compression, condensation, expansion, and evaporation processes, respectively. The processes 5′–5 and 3′–3 represent the sub-cooling of the refrigerant in the condenser and the superheating in the evaporator, respectively. Process 3–4s corresponds to the ideal compression process [31,32], while process 3–4 represents the actual compression process in the compressor. States 3′, 4′ and 5′ are the saturation states. The thermodynamic analysis was carried out using the Engineering Equation Solver Pro V10 (EES) software [38].

## 3. Energy, Exergy and Environment Impact Analysis

This section provides a detailed presentation of the set of equations for evaluating the energy and exergy efficiency of the system, its components, and the system and vapor compression cycle $COP$. The description of thermodynamic processes with the corresponding points has been presented in the previous section and in Figure 3. The equations presented in this section are derived based on the assumptions that the compression process is isentropic, the evaporation and condensation processes are isobaric, and the expansion process is isenthalpic [2–4]. Energy and exergy analyses are based on the 1st and 2nd laws of thermodynamics [31,32].

### 3.1. Evaporator

The heat source for GSHP is the ground, which is extracted using BHE. The temperature of the HTF at point 1 ($T_1$) was assumed to be equal to the source temperature ($T_{source}$). The approach ($T_{EAP}$) and superheating ($T_{SH}$) temperatures in the evaporator area were taken from experimental observations. It follows that the temperature in the compressor suction line ($T_3$) is given by Equation (1) [31,32]:

$$T_3 = T_1 - T_{EAP} = T_{source} - T_{EAP} = T_{evap} + T_{SH} \qquad (1)$$

In this case, the HTF temperature at point 2 is determined by Equation (2) [2,3]:

$$T_2 = T_1 - \frac{\dot{Q}_{evap}}{c_{p,1} \times \dot{m}_{pump-1}}. \qquad (2)$$

The amount of heat received by the evaporator from the low-grade ground heat source is:

$$\dot{Q}_{evap} = \dot{m}_{ref} \times (h_3 - h_6). \qquad (3)$$

### 3.2. Compressor

According to the definition of isentropic efficiency of the compressor, the refrigerant enthalpy after compression (state 4) can be calculated as follows:

$$h_4 = h_3 + \frac{h_{4s} - h_3}{\eta_{isen}}. \qquad (4)$$

In this case, the isentropic efficiency is determined by Equation (5) [39]:

$$\eta_{isen} = \frac{T_{evap} + 273.15}{T_{cond} + 273.15} + 0.0025 \times T_{evap}, \qquad (5)$$

The power and electricity consumption of the compressor can be determined as:

$$\dot{W}_{comp} = \dot{m}_{ref} \times (h_4 - h_3), \qquad (6)$$

$$\dot{W}_{elec,comp} = \frac{\dot{W}_{comp}}{\eta'_{comp}}, \qquad (7)$$

$$\eta'_{comp} = \frac{\eta_{comp}}{\eta_{isen}}, \qquad (8)$$

where $\eta'_{comp}$ is combination of mechanical, transmission and electrical motor efficiencies [31,32], and $\eta_{comp}$ is the ratio of ideal isentropic work to real electricity consumption. In this paper, coefficient $\eta_{comp}$ (in %) is taken as a function of $T_{evap}$ and $T_{cond}$ according to the compressor's manufacturer (Copeland ZR28K3E-PFJ) Equation (9):

$$\eta_{\text{comp}} = c_1 \times T_{\text{evap}} + c_2 \times T_{\text{cond}} + c_3 \times T_{\text{evap}}^2 + c_4 \times T_{\text{evap}} \times T_{\text{cond}} + c_5 \times T_{\text{cond}}^2 \\ + c_6 \times T_{\text{evap}}^3 + c_7 \times T_{\text{evap}}^2 \times T_{\text{cond}} + c_8 \times T_{\text{evap}} \cdot T_{\text{cond}}^2 \\ + c_9 \times T_{\text{cond}}^3 + c_{10}, \quad (9)$$

with $c_1 = -2.4341$, $c_2 = 1.367$, $c_3 = -0.062$, $c_4 = 0.103$, $c_5 = -0.0173$, $c_6 = -3.88 \times 10^{-5}$, $c_7 = 7.89 \times 10^{-4}$, $c_8 = -7.51 \times 10^{-4}$, $c_9 = 1.85 \times 10^{-6}$, $c_{10} = 33.6476$.

*3.3. Condenser*

The heat sink for GSHP is the water storage tank. The temperature of the HTF at point 7 ($T_7$) was taken equal to the sink temperature ($T_{\text{sink}}$). The approach ($T_{\text{CAP}}$) and sub-cooling ($T_{\text{SC}}$) temperatures in the condenser area were given by experimental observations. Thus, condenser outlet temperature can be determined by Equation (10) [31,32]:

$$T_5 = T_7 + T_{\text{CAP}} = T_{\text{sink}} + T_{\text{CAP}} = T_{\text{cond}} - T_{\text{SC}} \quad (10)$$

In this case, the temperature of the HTF at point 8 is determined by Equation (11) [2,3]:

$$T_8 = T_7 + \frac{\dot{Q}_{\text{cond}}}{c_{\text{p},7} \times \dot{m}_{\text{pump}-2}}. \quad (11)$$

The amount of useful heat from the condenser is determined by Equation (12):

$$\dot{Q}_{\text{cond}} = \dot{m}_{\text{ref}} \times (h_4 - h_5). \quad (12)$$

*3.4. Expansion Device*

Since the expansion process is assumed to take place at a constant enthalpy, the enthalpy at the exit of the expansion valve is equal to that at the inlet, i.e., $h_6 = h_5$.

*3.5. System Parameters*

The refrigerant mass flow rate is obtained from the condenser heating capacity by Equation (13):

$$\dot{m}_{\text{ref}} = \frac{\dot{Q}_{\text{cond}}}{h_4 - h_5}. \quad (13)$$

The overall electrical power input for running the system is the sum of the power consumptions of the compressor and the two circulating pumps:

$$\dot{W}_{\text{elec-total}} = \dot{W}_{\text{elec,comp}} + \dot{W}_{\text{elec,pump}-1} + \dot{W}_{\text{elec,pump}-2}. \quad (14)$$

The volumetric refrigerating capacity at the outlet of the evaporator is calculated as:

$$VRC = \rho_3 \times (h_3 - h_6). \quad (15)$$

The $COP$ of the vapor compression heat pump cycle (process 3-4-5-6) and the $COP$ of the GSHP system are given by Equations (16) and (17), respectively:

$$COP_{\text{cycle}} = \frac{\dot{Q}_{\text{cond}}}{\dot{W}_{\text{elec,comp}}}, \quad (16)$$

$$COP_{\text{system}} = \frac{\dot{Q}_{\text{cond}}}{\dot{W}_{\text{elec-total}}}. \quad (17)$$

The $COP$ of the reverse Carnot cycle is expressed in terms of the source and sink temperatures as:

$$COP_{\text{carnot}} = \frac{T_{\text{sink}}}{T_{\text{sink}} - T_{\text{source}}} = \frac{T_7}{T_7 - T_1}. \quad (18)$$

*3.6. Exergy Destruction*

Exergy analysis provides a method for estimating the maximum work extracted from a substance relative to its environment as a reference state (the "dead state") [2,3]. Exergy analysis is widely used as a method for the thermodynamic optimization of thermal energy systems, which helps determine the inefficient components of GSHP [31,32].

Exergy destruction of vapor compression heat pump cycle components is determined by Equation (19), which is derived from the steady-state form of the exergy rate balance:

$$\dot{E}^{\text{dest}} = \sum \dot{E}^{\text{in}} - \sum \dot{E}^{\text{out}} - \dot{W} + \dot{E}_Q. \tag{19}$$

The first and second terms of the right-hand-side of Equation (19) are the flow exergy at the inlet and outlet of a specific system component; the third term is related to the work done by the component, and the fourth term corresponds to the heat exchange between environment and the component.

Specific flow exergy and flow exergy at the $i$th point are determined based on enthalpy and entropy of this point and the dead state:

$$e_i = (h_i - h_0) - T_0 \times (s_i - s_0), \tag{20}$$

$$\dot{E}_i = \dot{m}_i \times e_i = \dot{m}_i \times [(h_i - h_0) - T_0 \times (s_i - s_0)]. \tag{21}$$

For each component of the heat pump, the exergy destruction is given by Equation (22) to Equation (25) for the evaporator, compressor, condenser, and expansion valve, respectively:

$$\dot{E}^{\text{dest}}_{\text{evap}} = (\dot{E}_1 + \dot{E}_6) - (\dot{E}_2 + \dot{E}_3), \tag{22}$$

$$\dot{E}^{\text{dest}}_{\text{comp}} = \dot{E}_3 + \dot{W}_{\text{elec,comp}} - \dot{E}_4, \tag{23}$$

$$\dot{E}^{\text{dest}}_{\text{cond}} = (\dot{E}_4 + \dot{E}_7) - (\dot{E}_5 + \dot{E}_8), \tag{24}$$

$$\dot{E}^{\text{dest}}_{\text{TEV}} = \dot{E}_5 - \dot{E}_6. \tag{25}$$

The overall cycle exergy destruction is determined by Equation (26), as the sum of $j$th components exergy destruction:

$$\dot{E}^{\text{dest}}_{\text{cycle}} = \sum \dot{E}^{\text{dest}}_j = \dot{E}^{\text{dest}}_{\text{comp}} + \dot{E}^{\text{dest}}_{\text{cond}} + \dot{E}^{\text{dest}}_{\text{TEV}} + \dot{E}^{\text{dest}}_{\text{evap}}. \tag{26}$$

The relative exergy destruction of the $j$th component is determined by Equation (27):

$$\text{RE}^{\text{dest}}_j = \frac{\dot{E}^{\text{dest}}_j}{\dot{E}^{\text{dest}}_{\text{cycle}}}. \tag{27}$$

*3.7. Exergy Efficiency*

The evaporator and condenser are brazed plate heat exchangers, in which heat is transferred with counter-current flows (see Table 3). The exergy efficiency of the evaporator and condenser can be determined by Equations (28) and (29), respectively:

$$\eta^{\text{ex}}_{\text{evap}} = \frac{\dot{E}_1 - \dot{E}_2}{\dot{E}_3 - \dot{E}_6}. \tag{28}$$

$$\eta^{\text{ex}}_{\text{cond}} = \frac{\dot{E}_8 - \dot{E}_7}{\dot{E}_4 - \dot{E}_5}. \tag{29}$$

The vapor compression process of the refrigerant in the compressor is assumed to be adiabatic. Therefore, its exergy efficiency can be determined by Equation (30):

$$\eta_{\text{comp}}^{\text{ex}} = \frac{\dot{E}_4 - \dot{E}_3}{\dot{W}_{\text{elec,comp}}}. \tag{30}$$

The expansion process of the refrigerant is assumed to be isenthalpic. The exergy efficiency of the expansion valve can be determined by the ratio of the exergy output to the exergy input [40]:

$$\eta_{\text{TEV}}^{\text{ex}} = \frac{\dot{E}_6}{\dot{E}_5}. \tag{31}$$

The exergy efficiencies of the heat pump cycle and GSHP system can be calculated by Equations (32) and (33) [2,3]:

$$\eta_{\text{cycle}}^{\text{ex}} = \frac{\dot{E}_8 - \dot{E}_7}{\dot{W}_{\text{elec,comp}}}, \tag{32}$$

$$\eta_{\text{system}}^{\text{ex}} = \frac{\dot{E}_8 - \dot{E}_7}{\dot{W}_{\text{elec-total}}}. \tag{33}$$

The second-law efficiency of the cycle is defined as the ratio of the cycle $COP$ to the $COP$ of the reverse Carnot cycle [31,32]:

$$\eta_{2\text{nd}} = \frac{COP_{\text{cycle}}}{COP_{\text{carnot}}}. \tag{34}$$

*3.8. Total Equivalent Warming Impact*

The total equivalent warming impact ($TEWI$) includes direct emission due to the system refrigerant leakage losses and indirect emission due to the heat pump electricity consumption. $TEWI$ is assessed by Equation (35) [41]:

$$TEWI = \underbrace{GWP \times m \times \left(\frac{L}{100} \times n + \left(1 - \frac{\alpha}{100}\right)\right)}_{\text{direct}} + \underbrace{E \times \beta \times n}_{\text{indirect}}, \tag{35}$$

where $GWP$ is the Global Warming Potential of the refrigerant (in $CO_2$ eq.), $m$ is the charge of the refrigerant (in kg), $L$ is the leakage rate per year (in %/year), $n$ is system operating life (in years); $\alpha$ is the recovery/recycling factor (in %), $E$ is the energy consumption per year (in kWh/year), $\beta$ is the indirect emission factor (in kg/kWh).

## 4. Results and Discussions

This section presents experimental observations, validation of theoretical models, and thermodynamic performance comparison of water-to-water heat pump using environment friendly refrigerants. In addition, the total equivalent warming impact of the system working with R134a, and its possible alternatives is discussed.

*4.1. Experimental Results*

The measurements were carried out using the measuring devices shown in Figure 2. As seen from Figure 2, temperature and pressure sensors are installed at the nodal points of the GSHP. In addition, in the hydraulic loops between BHE and the evaporator, as well as between the condenser and the water storage tank, liquid flow meters are installed. The observations are depicted in Figures 4–8 for the water heating in the storage tank.

*Energies* **2022**, *15*, 8336　　　　　　　　　　　　　　　　　　　　　　　　　　　11 of 25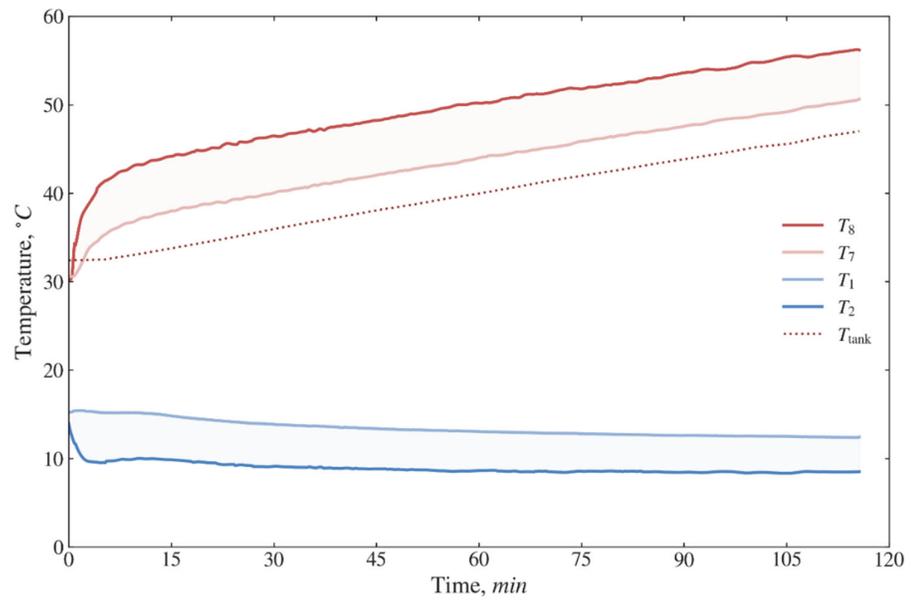

**Figure 4.** Temperatures measured for the heat transfer fluid vs time.

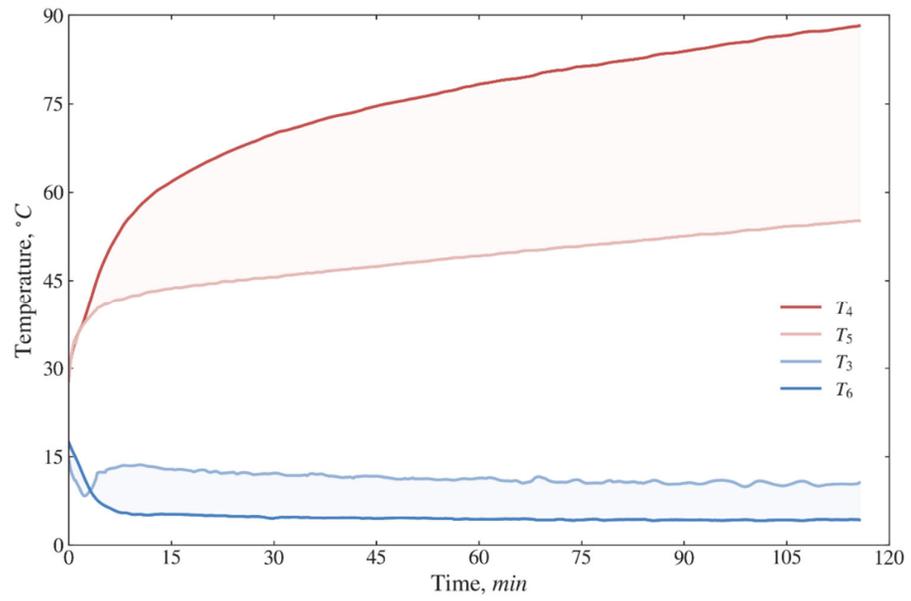

**Figure 5.** Temperatures measured for refrigerant R134a vs time.

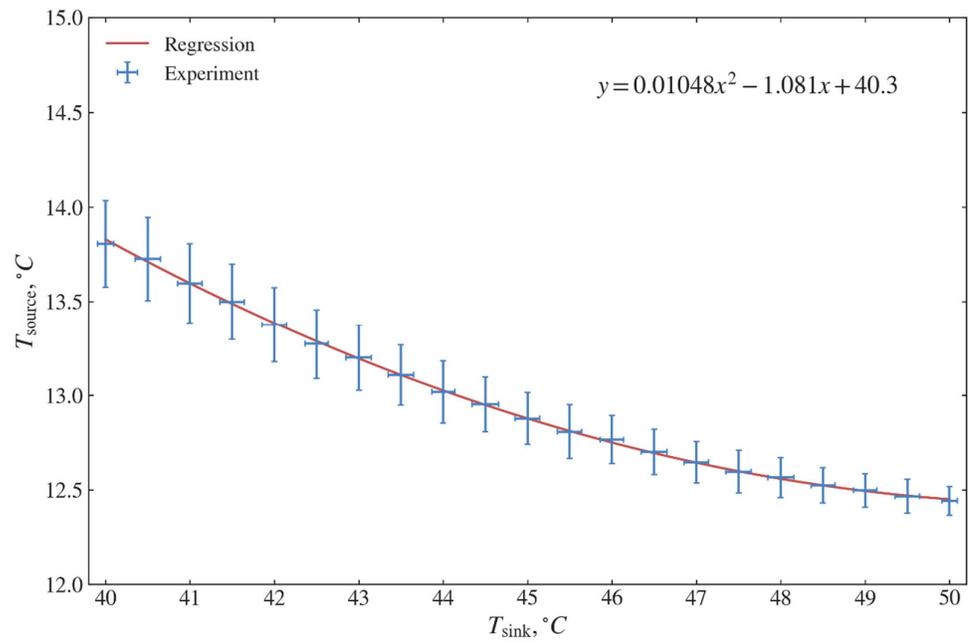

**Figure 6.** Variation of the temperature of the heat source with the temperature of the heat sink.

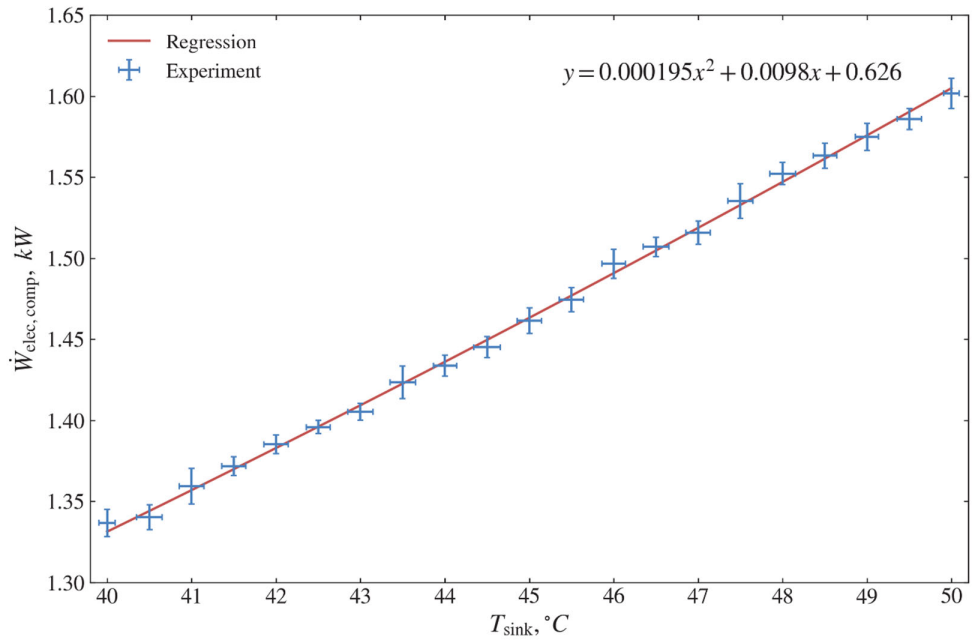

**Figure 7.** Variation of compressor electricity consumption measured with the temperature of the heat sink.

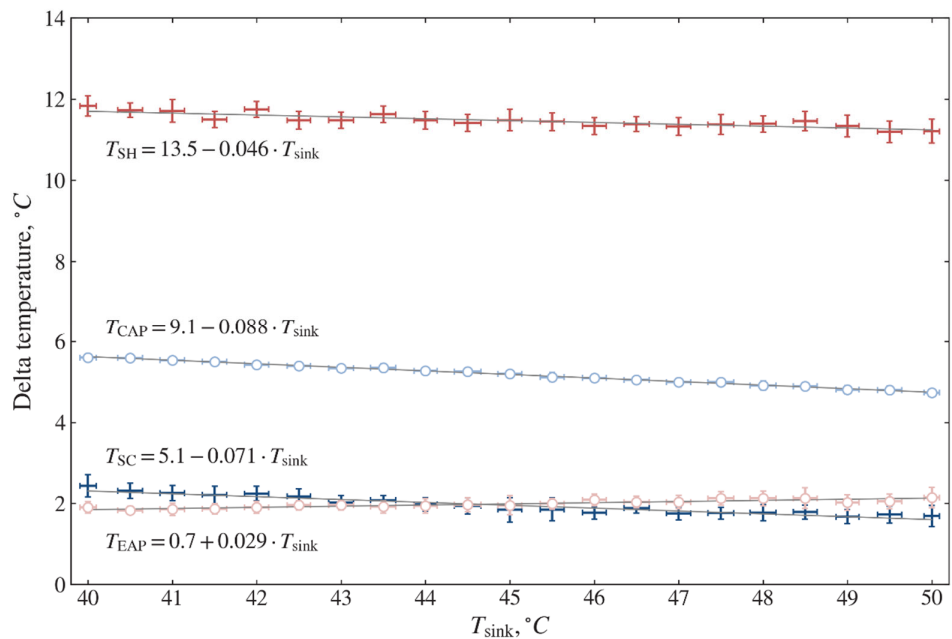

**Figure 8.** Variations of estimated delta temperatures based on measured temperature with the temperature of the heat sink.

Figures 4 and 5 show the measurement data of a single experiment for the temperatures of the HTF and refrigerant R134a ( see the nodal points shown in Figure 2). Figure 4 shows the results on changes in the water temperature in the tank. For two hours of operation of the GSHP, after establishing the vapor compression cycle operation (≈ 10 min) the temperature range for $T_1$ is 15.2–12.4 °C, $T_2$ is 10–8.5 °C, $T_7$ is 37–50.6 °C, $T_8$ is 43.2–56.2 °C, while the temperature of the water in the tank rises from 33 °C to 47 °C. As shown in Figure 5, the temperatures of refrigerant R134a at the nodal points of the vapor compression heat pump cycle vary in the following ranges: $T_3$ = 13.6–10.6 °C, $T_4$ = 57–88 °C, $T_5$ = 42.5–55.2 °C, $T_6$ = 5.2–4.3 °C.

To describe the operation of GSHP, as indicated in the diagram of Figure 2 and in the thermodynamic model (1)–(35), the temperature of the HTF $T_1$ was taken as the source temperature $T_{source}$, while the temperature $T_7$ was taken as the sink temperature $T_{sink}$. All measurements and calculations were obtained for $T_{sink}$ values in the range of 40–50 °C. The results of several experiments on $T_{source}$ and $\dot{W}_{elec,comp}$ are presented in Figures 6 and 7, respectively.

As shown in Figure 6, the source temperature changes from 12.45 °C to 13.83 °C. The ground type of the local terrain is mainly composed of gravel and till/loam. For the specified type of ground, the operation of the constructed GSHP gives the specified regression dependence between two temperatures, which is then used for validation (see Table 4). A similar relationship for measuring the compressor electricity consumption $\dot{W}_{elec,comp}$ is shown in Figure 7. $\dot{W}_{elec,comp}$ is in the range of 1.33–1.61 kW.

Table 4. Model operation conditions.

| Parameter | Designation | Unit | Value * |
|---|---|---|---|
| Heat sink temperature | $T_{sink}$ | °C | 40–50 {50} |
| Heat source temperature | $T_{source}$ | °C | $T_{source} = 0.01048 \times T_{sink}^2 - 1.081 \times T_{sink} + 40.3$ {13} |
| Compressor electricity consumption | $\dot{W}_{elec,comp}$ | kW | $\dot{W}_{elec,comp} = 0.000195 \times T_{sink}^2 + 0.0098 \times T_{sink} + 0.626$ |
| Compressor efficiency | $\eta_{comp}$ | % | Equation (9) |
| Superheating temperature | $T_{SH}$ | °C | $T_{SH} = 13.5 - 0.046 \times T_{sink}$ {8} |
| Sub-cooling temperature | $T_{SC}$ | °C | $T_{SC} = 5.1 - 0.071 \times T_{sink}$ {2} |
| Evaporator approach temperature | $T_{EAP}$ | °C | $T_{EAP} = 0.7 + 0.029 \times T_{sink}$ {2} |
| Condenser approach temperature | $T_{CAP}$ | °C | $T_{CAP} = 9.1 - 0.088 \times T_{sink}$ {2} |
| Dead state temperature | $T_0$ | °C | 20 |
| Dead state pressure | $P_0$ | kPa | 101.325 |
| Pump1/pump2 volume flow rate | $\dot{V}_{pump-1}$ | m³/h | 0.595 |
| | $\dot{V}_{pump-2}$ | | 0.543 |
| Pump1/pump2 power consumption | $\dot{W}_{elec,pump-1}$ | W | 60 |
| | $\dot{W}_{elec,pump-2}$ | | 40 |

* The values in brackets { } are used for the calculations in Section 4.3.

Furthermore, to compare the measured and calculated results, the following parameters such as superheating $T_{SH}$, sub-cooling $T_{SC}$, evaporator approach temperature $T_{EAP}$, and condenser approach temperature $T_{CAP}$ are used and encountered in Equations (1) and (10). Based on the superheating and sub-cooling results, the operation of the vapor compression heat pump system, the amount of refrigerant mass charge in the system, and the thermostatic expansion valve, adjustments are ascertained. Figure 8 shows the calculated values of $T_{SH}$, $T_{SC}$, $T_{EAP}$, and $T_{CAP}$ based on temperature measurements. For these parameters, linear regression dependencies are found. The values are in the ranges of $T_{SH}$ = 11.2–11.7 °C, $T_{SC}$ = 1.6–2.3 °C, $T_{EAP}$ = 1.8–2.1 °C, $T_{CAP}$ = 4.7–5.6 °C. These values are used for validation (see Table 4).

*4.2. Thermodynamic Model Validation*

This subsection contains the comparison of experimental and theoretical results on temperatures in the heat pump cycle, evaporator, and condenser heat fluxes, $COPs$ of the cycle and system. For the thermodynamic calculations based on Equations (1)–(18), the parameters indicated in Table 4 are taken. The theoretically predicted thermodynamic performance results are compared with experimental results for validation.

Figure 9 illustrates the temperature at typical nodal points as the function of the sink temperature $T_{sink}$. The experimental and calculated values of the refrigerant temperatures at the compressor suction line $T_3$, as well as at the condenser outlet $T_5$, are in very good agreement. The temperature difference at the evaporator inlet $T_6$ is 3.16–4.01 °C with an average value 3.66 °C. The difference at the compressor discharge line $T_4$ is 0.06–6.12 °C. This difference between the experiment and the model is because the model does not consider the presence of a four-way reverse valve.

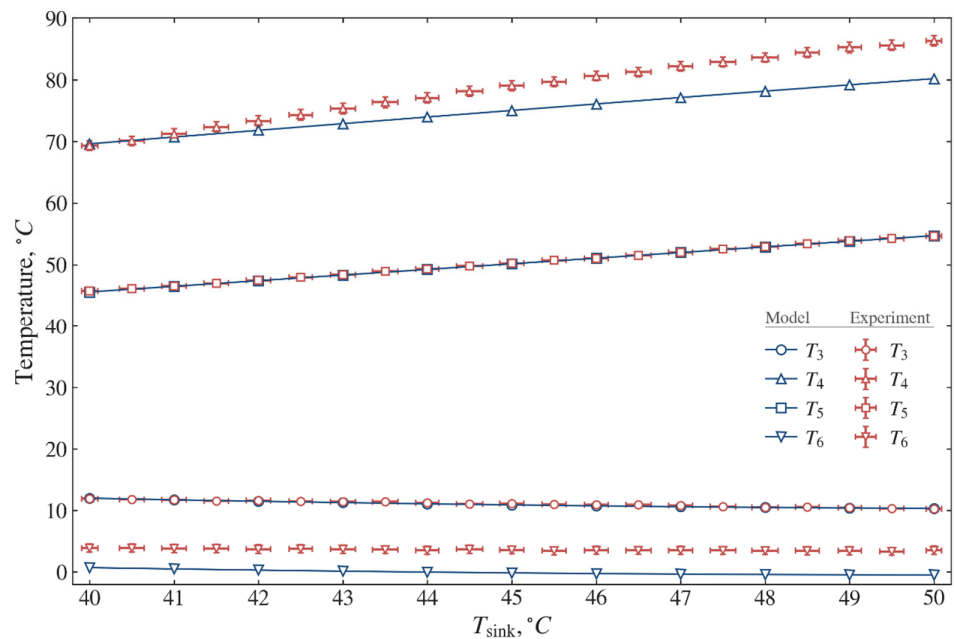

**Figure 9.** Comparisons of predicted and measured temperatures with the temperature of the heat sink.

Figure 10 shows comparisons of the predicted and measured heat capacities for the evaporator and condenser. The measured amount of useful heat $\dot{Q}_{cond}$ provided by the condenser is in the range of 3.59–3.96 kW, while the simulated values are in the range of 3.88–4.26 kW. The average systematic difference is 6.2%, which is equivalent to 251 W. Such a systematic difference is because the presented thermodynamic model does not consider pressure drops and heat losses. The experimental value of received heat by evaporator $\dot{Q}_{evap}$ is in the range of 2.73–3.31 kW, while the simulated values are in the range of 2.86–3.37 kW. The average systematic difference is 2.6%, which is equivalent to 79.6 W.

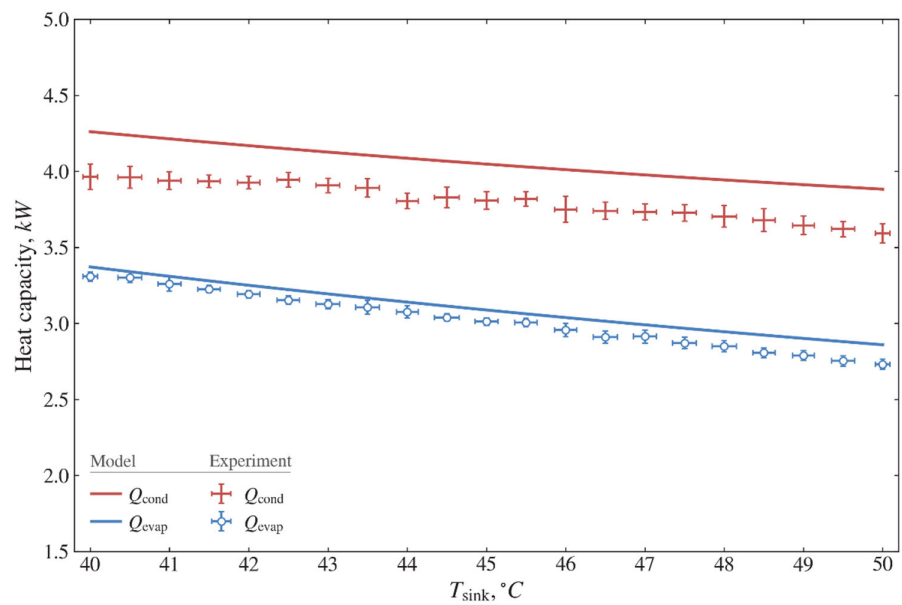

**Figure 10.** Comparisons of predicted and measured heat capacities with the temperature of the heat sink for the condenser and evaporator.

Figure 11 shows the comparisons of the predicted and measured $COPs$ of the cycle and system with the temperature of the heat sink. The experimental values of $COP_{cycle}$ vary in the range of 2.24–2.97, while theoretical values calculated using Equation (16) vary in the range of 2.42–3.20. The average difference between the experimental and theoretical values is 6.2%. The experimental values of $COP_{system}$ vary in the range of 2.11–2.76, while the theoretical values calculated using Equation (17) vary in the range of 2.28–2.98. The minimum, maximum and average differences between the experimental and theoretical values are 4.9%, 7.5% and 6.2%, respectively. Therefore, the thermodynamic model developed in the present work is validated and can be used to further analyze the thermal performance of the GSHP system.

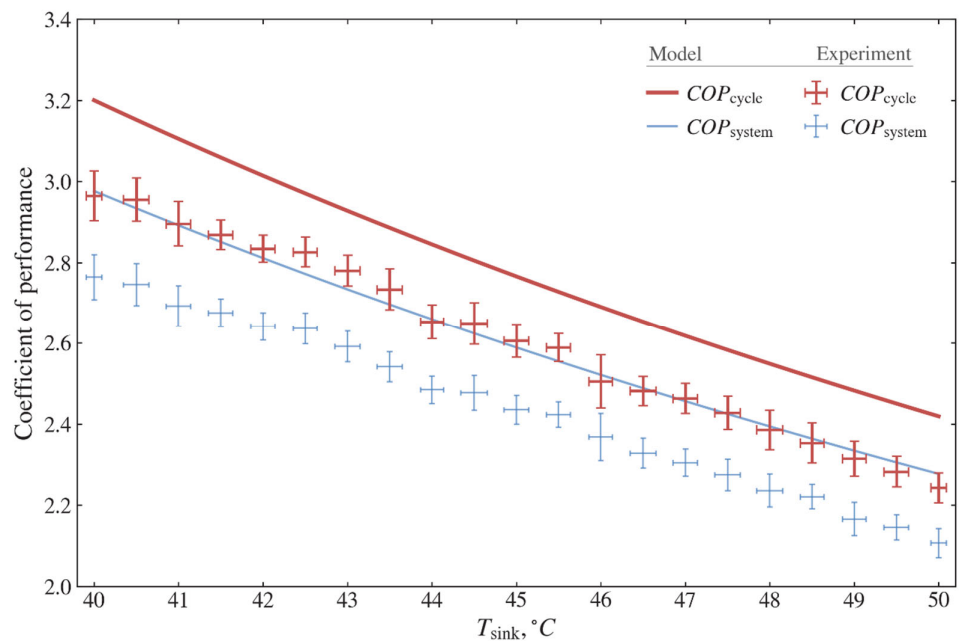

**Figure 11.** Comparisons of predicted and measured $COPs$ of the cycle and system with the temperature of the heat sink.

*4.3. Thermodynamic and Environmental Calculations*

Using the theoretical model, the thermodynamic performance of the heat pump is simulated for the environment friendly refrigerant options. The principle behind selecting alternative refrigerants is affordability, similar thermodynamic and thermo-physical properties, low $GWP$ values, and non- or mild flammability. Table 5 shows the thermodynamic properties of selected refrigerants in this work.

**Table 5.** Thermodynamic properties of the selected refrigerants.

| Property | Unit | R134a | R152a | R450A | R513A | R1234yf | R1234ze(E) |
|---|---|---|---|---|---|---|---|
| $GWP$ for 100 years | $CO_2$ eq. | 1430 | 124 | 547 | 573 | 4 | 7 |
| Molar mass | g/mol | 102.03 | 66.05 | 108.6 | 108.4 | 114.0 | 114.0 |
| Normal boiling point | °C | −26.1 | −24.7 | −23.1 | −29.2 | −29.4 | −19.0 |
| Critical pressure | MPa | 4.06 | 4.50 | 4.01 | 3.77 | 3.38 | 3.64 |
| Critical temperature | °C | 101.1 | 113.15 | 75.1 | 96.5 | 94.7 | 109.4 |
| UFL | vol. % | - | 16.9 | - | - | 12.3 | 11.3 |
| LFL | vol. % | - | 3.9 | - | - | 6.2 | 5.7 |
| Auto-ignition temperature | °C | - | 455 | - | - | 405 | 368 |

The thermodynamic analysis is performed using the model of Equations (1)–(35) and the operating conditions given in Table 4. However, instead of the given input power consumption, a value of $\dot{Q}_{cond}$ = 5 kW is fixed as an input condition. Figure 12 shows the relationship between $\dot{m}_{ref}$ required of the system and the specific latent heat of vaporization of refrigerants. The highest value of the specific latent heat of vaporization corresponds to refrigerant R152a with the lowest value of $\dot{m}_{ref}$ while the lowest value of the specific latent heat refers to refrigerant R1234yf with the highest value of $\dot{m}_{ref}$. As seen from Figure 12, the refrigerant mass flow rate required with the increase in the specific latent heat of vaporization of refrigerants.

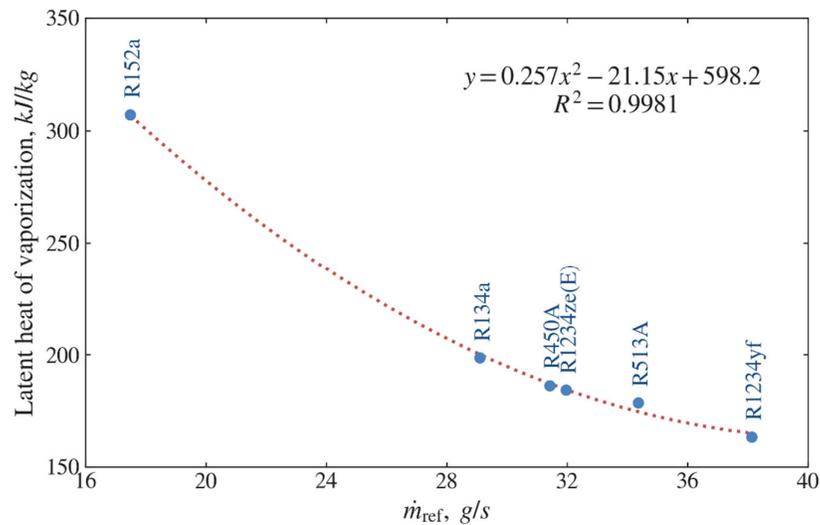

**Figure 12.** Refrigerant mass flow rate required for a given heat capacity vs latent heat of vaporization of refrigerants.

Figure 13 shows the comparison of the required refrigerant mass flow rate and volumetric refrigeration capacity ($VRC$) for a given heat capacity for refrigerant R134a and environment-friendly alternative refrigerants. For the same heating capacity as refrigerant R134a, the low refrigerant mass flow rate ($\dot{m}_{ref}$) with high $VRC$ results in a smaller size compressor. In this regard, refrigerant R152a is seen to be the best option. For R152a, the refrigerant mass flow rate ($\dot{m}_{ref}$) required is estimated to be 17.5 g/s and $VRC$ 1949 kJ/m³. The mass flow rates required for R450A and R1234ze(E) are close to that of R134a, with the values of 31.4 g/s and 32.0 g/s, respectively. However, refrigerant R450A and R1234ze(E) have lower values of $VRC$ than that of R134a, which requires a large volume compressor. In addition, R513A and R1234yf have closer values of $VRC$ to that of R134a at higher mass flow rates.

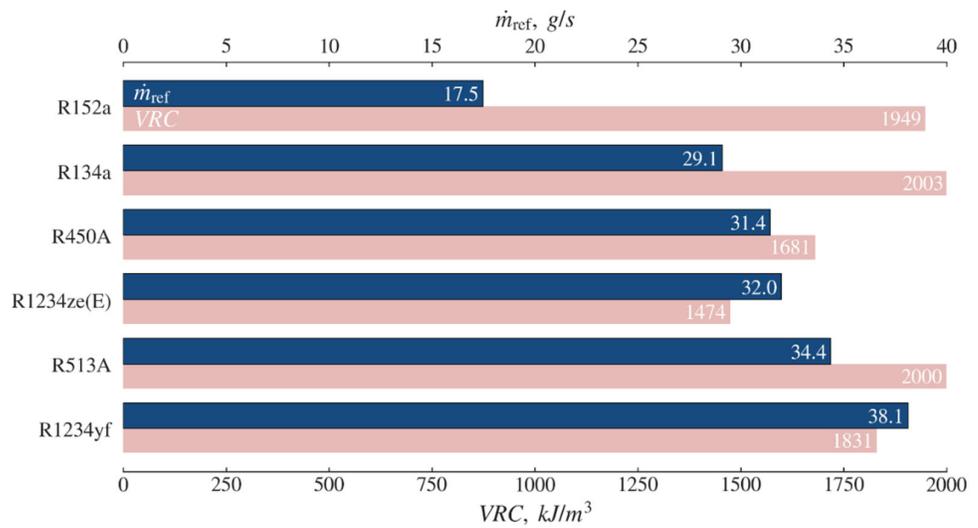

**Figure 13.** Refrigerant mass flow rate and volumetric refrigeration capacity required for a given heat capacity for different refrigerants.

Figures 14 and 15 show the comparison of the exergy destructions and relative exergy destructions of the system components and cycle for different refrigerants, respectively. Among the selected refrigerants, the minimum system exergy destruction $\dot{E}^{dest}_{system}$ is observed for R152a with a value of 965 W and the maximum exergy destruction of 1119 W is observed for R1234yf. Among components, the maximum exergy destruction is observed in the compressor, followed by the expansion valve, evaporator, and condenser.

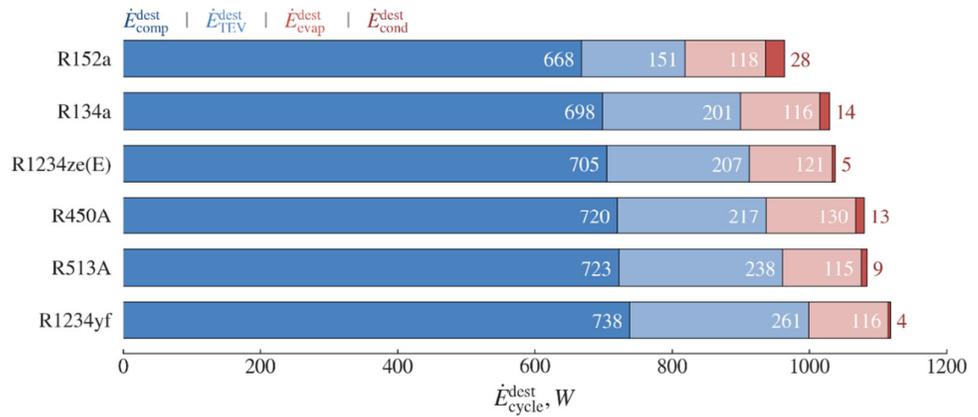

**Figure 14.** Exergy destructions of the components and cycle for different refrigerants.

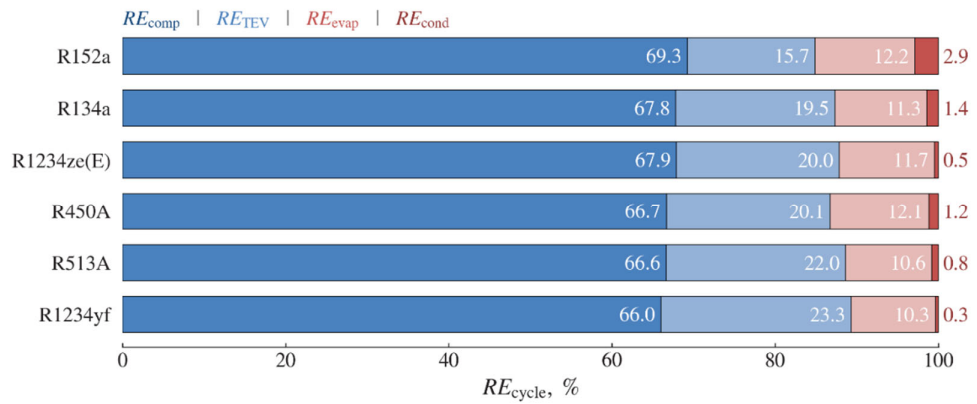

**Figure 15.** Relative exergy destructions of the components and cycle for different refrigerants.

As is seen from Figure 15, the expansion valve can apparently affect the growth of system exergy destruction for these refrigerants. The share of the expansion device $RE_{TEV}^{dest}$ increases from 16% to 23%, while the relative exergy destructions of other components vary slightly by 1–3%. Referring to Table 4, since $T_{source}$ and HTF flow rate $\dot{V}_{pump-1}$ are fixed for the evaporator, for these refrigerants, there is no significant difference in $\dot{E}_{evap}^{dest}$. Although $T_{sink}$ and $\dot{V}_{pump-2}$ for the condenser are fixed, $\dot{E}_{cond}^{dest}$ is different. This is due to the different discharge temperature $T_4$ depending on the thermodynamic properties of the refrigerant.

Figure 16 shows the effect of discharge temperature $T_4$ on the condenser exergy efficiency $\eta_{cond}^{ex}$. The highest value of $\eta_{cond}^{ex}$ corresponds to refrigerant R1234yf with the lowest value of $T_4$ whereas the lowest value of $\eta_{cond}^{ex}$ refers to refrigerant R152a with the highest value of $T_4$. It is found that a linear relationship between $T_4$ and $\eta_{cond}^{ex}$. Thus, this dependence shows that with an increase in the discharge temperature, the condenser efficiency decreases in this operation mode. For these refrigerants, the condenser exergy efficiency varies between 94.9% and 99.3%. The exergy efficiencies of other components are given in Table 6. The exergy efficiencies of the evaporator, compressor and expansion valve are in the range of 50.4–53.6%, 58.3–58.8% and 79.9–85.7%, respectively.

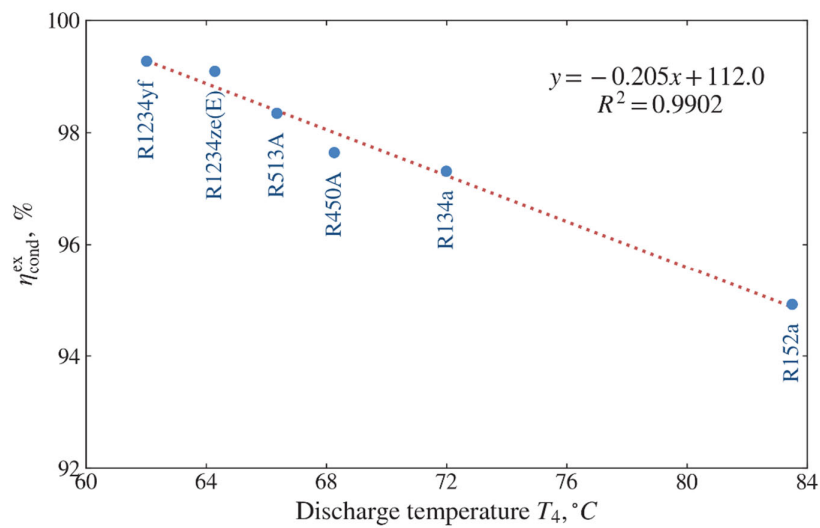

**Figure 16.** Condenser exergy efficiency vs discharge temperature for different refrigerants.

The second law $\eta_{2nd}$, cycle $\eta_{cycle}^{ex}$, and overall system $\eta_{system}^{ex}$ efficiencies are calculated using Equations (32)–(34). Figure 17 shows the evaluating results of these efficiencies. These efficiencies for all the refrigerants appear a similar trend in the ranges of 32.3–35.3%, 29.4–32.1%, and 27.8–30.2%, respectively. As is seen from Table 6, the maximal values of the efficiencies correspond to refrigerant R152a, and the minimal values to refrigerant R1234yf. Among the selected refrigerants, only R152a is superior to R134a while R1234ze(E) has the closest values to those of R134a. The values of $COP_{cycle}$ and $COP_{system}$ vary similarly.

Table 6. The efficiencies of heat pump cycle and components.

| Refrigerant | $COP_{cycle}$ | $COP_{system}$ | $\eta_{2nd}$, % | $\eta_{cycle}^{ex}$, % | $\eta_{system}^{ex}$, % | $\eta_{evap}^{ex}$, % | $\eta_{comp}^{ex}$, % | $\eta_{cond}^{ex}$, % | $\eta_{TEV}^{ex}$, % |
|---|---|---|---|---|---|---|---|---|---|
| R152a | 3.09 | 2.91 | 35.3 | 32.1 | 30.2 | 53.6 | 58.8 | 94.9 | 85.7 |
| R134a | 2.97 | 2.80 | 34.0 | 30.9 | 29.1 | 53.6 | 58.5 | 97.3 | 83.5 |
| R1234ze(E) | 2.96 | 2.79 | 33.8 | 30.7 | 29.0 | 52.5 | 58.3 | 99.1 | 79.7 |
| R450A | 2.89 | 2.73 | 33.1 | 30.0 | 28.4 | 50.4 | 58.4 | 97.6 | 81.3 |
| R513A | 2.88 | 2.72 | 32.9 | 29.9 | 28.3 | 53.6 | 58.4 | 98.3 | 83.1 |
| R1234yf | 2.83 | 2.67 | 32.3 | 29.4 | 27.8 | 53.2 | 58.3 | 99.3 | 82.1 |

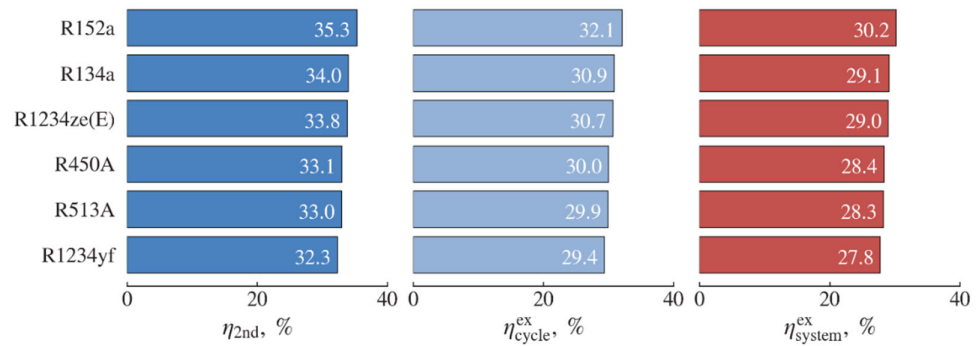

**Figure 17.** Second law, cycle, and system efficiencies for different refrigerants.

Figure 18 shows the total equivalent warming impact ($TEWI$) of a water-to-water heat pump using R134a and environment-friendly alternative refrigerants. The $TEWI$ includes both direct and indirect emissions depending on the $GWP$ of the refrigerant. R134a appears the maximum direct emission due to its high $GWP$ than other refrigerants, followed by R513A, R450A, and R152a having the direct impacts of 1.14, 0.84, and 0.12 tons of $CO_2$ equivalent, respectively. The lowest value of direct emissions is observed for HFO refrigerants like R1234yf and R1234ze(E) with a value of 0.01. The lowest value of indirect emissions is observed for R152a with 31.64-ton eq. $CO_2$ while the highest value of 34.32-ton eq. $CO_2$ is associated with R1234yf due to its high compressor power consumption. The carbon-dioxide emission factor ($\beta$) is assumed to be 0.8 kg of $CO_2$ per kWh for coal-based electric power generation system [41].

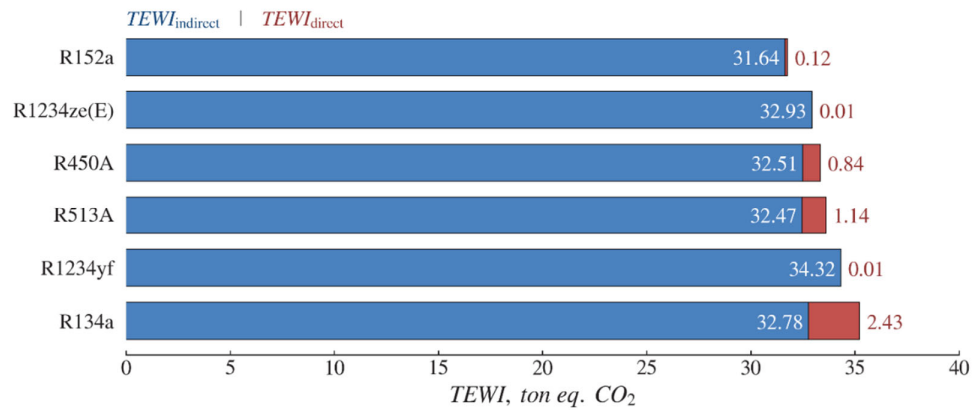

**Figure 18.** Total equivalent warming impact for different refrigerants.

Figure 19 shows normalized system parameter dependences like $COP_{cycle}$ vs. $TEWI$ and $\eta^{ex}_{cycle}$ vs. $\dot{E}^{dest}_{cycle}$ for different refrigerants. R134a is used as the base refrigerant and others refrigerants are considered as the low $GWP$ alternatives. HFCs such as R134a will be phased out soon according to the Paris protocol 2016 (Kigali agreement). The second quadrant is the best choice for both cases in Figure 19 (a) and (b). In this respect, a normalized $TEWI$ lower than one, and a $COP_{cycle}$ greater than one is preferable. Similarly, a normalized $\eta^{ex}_{cycle}$ greater than one, and a $\dot{E}^{dest}_{cycle}$ lower than one is preferable.

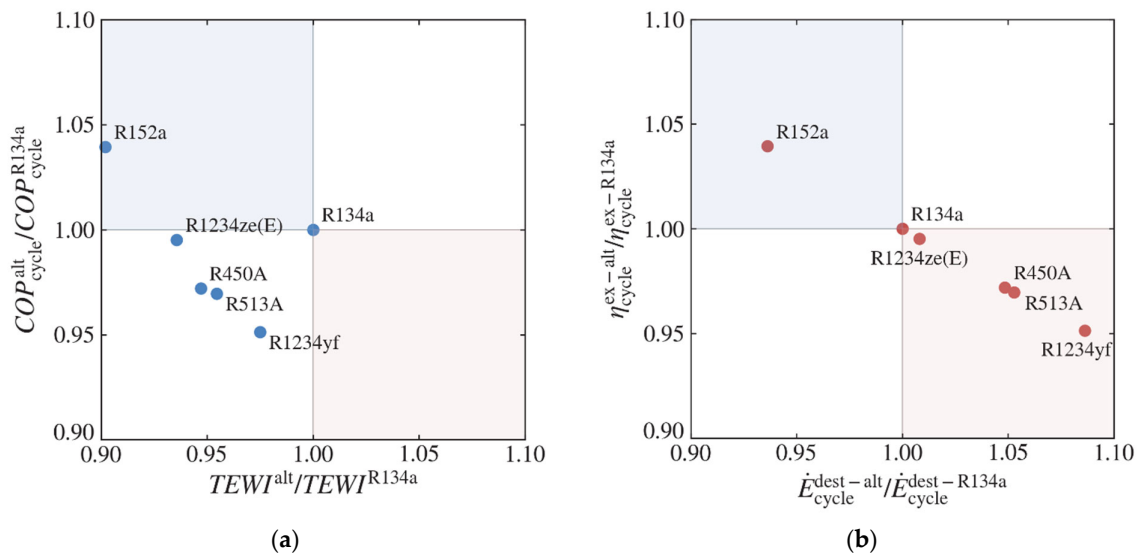

**Figure 19.** Normalized cycle parameter dependences for different refrigerants: (**a**) Coefficient of performance vs. total equivalent warming impact; (**b**) Cycle efficiency vs. exergy destruction.

As shown in Figure 19, R152a is a good choice with higher $COP_{cycle}$ and $\eta^{ex}_{cycle}$ than R134a; moreover, its values of $TEWI$ and $\dot{E}^{dest}_{cycle}$ are the lowest. The next refrigerant with lower $TEWI$ is R1234ze(E), while the remaining parameters are close to those of R134a. The lowest efficiency parameters are observed for R1234yf with a lower $TEWI$ than that of R134a, but with the highest exergy destruction. Refrigerant R450A and R513A display similar properties, with efficiencies lower than that of R134a despite their low $TEWI$.

Refrigerant R152a is seen to be a good alternative to replace R134a. However, due to its flammability, the mass charge of R152a is limited [40]. The HFO refrigerants such as R1234yf and R1234ze(E) are considered to be the key substitutes for the HFC refrigerants.

Even though these HFOs have low global warming potential, the thermo-physical properties are found to be lower than those of R134a, which results in poor thermodynamic performance. Moreover, fast decomposition of HFO refrigerants in the atmosphere results in formation of trifluoroacetate (TFA) [42]. For example, the level of TFA for R1234yf is almost five times higher than that for R134a [42]. The release of TFA in the atmosphere leads to the pollution of aquatic systems [2,3]. Therefore, HFC/HFO blends such as R450A and R513A are being considered as alternative refrigerants. Refrigerant R450A and R513A are non-flammable with a $GWP$ 60% lower than that of R134a. For these reasons, R450A and R513A can be recommended as good alternatives.

## 5. Conclusions

A prototype of a ground source heat pump heating system has been developed and installed in the Koksai Mosque, in the Almaty region of Kazakhstan. Experiments have been performed for refrigerant R134a. The energy and exergy efficiencies of the ground source heat pump heating system operating under the weather conditions of Kazakhstan have been analyzed. The thermodynamic model was validated by comparison with experimental data. Based on the verified model, calculation of energy and exergy efficiencies for environmentally friendly refrigerants were carried out as alternatives to R134a. The environmental impact of the heating system using different refrigerants has been estimated using the $TEWI$ parameter. The following major conclusions can be drawn from this study:

- The experimental values of $COP_{cycle}$ vary in the range of 2.24–2.97, while $COP_{system}$ vary in the range of 2.11–2.76.
- The predicted and experimental results are found to be in good agreement within 6.2%.
- The thermodynamic model describes the energy and exergy efficiencies of the system with sufficient accuracy.
- As input data for model calculations, we have used measurements from the borehole heat exchanger, that takes into account the local ground conditions.
- The maximal exergy destruction is attributed to the compressor, followed by the expansion valve, evaporator, and condenser.
- The exergy efficiency of the condenser depends on the compressor discharge temperature for different refrigerants, and it is maximal in the case of R1234yf.
- The R152a is found to be a good alternative to R134a in terms of $COP_{cycle}$, $COP_{system}$, and $TEWI$ with estimated values of 3.09, 2.91, and 31.76-ton eq. $CO_2$, respectively.
- Refrigerant R450A and R513A are found to be safer alternatives to R134a with slightly (of order 3%) lower $COP_{cycle}$, $COP_{system}$, and $TEWI$.
- Refrigerant R1234yf and R1234ze are identified as interim replacements to R134a due to the formation of trifluoroacetic acid.

**Author Contributions:** Conceptualization, Y.Y., A.T. (Amankeldy Toleukhanov), M.M., O.B., M.F., H.S.W. and Y.B.; methodology, O.B., M.F. and Y.B.; software, Y.Y.; validation, Y.Y., T.A., A.A., A.S. and A.T. (Amankeldy Toleukhanov); formal analysis, A.A. and A.S.; investigation, Y.Y., H.S.W. and Y.B.; resources, A.T. (Amankeldy Toleukhanov), O.B., H.S.W., A.T. (Alexandr Tsoy) and Y.B.; data curation, Y.Y., A.A. and A.S.; writing—original draft preparation, Y.B.; writing—review & editing, M.M., O.B., M.F. and H.S.W.; visualization, Y.Y. and T.A.; supervision, M.M., M.F. and H.S.W.; project administration, A.T. (Amankeldy Toleukhanov) and Y.B. All authors have read and agreed to the published version of the manuscript.

**Funding:** This research is funded by the Committee of Science of the Ministry of Science and Higher Education of the Republic of Kazakhstan, Grant No. AP08857319 "Study of heat transfer enhancement mechanisms of vertical type borehole heat exchanger to ensure high heat pump performance".

**Data Availability Statement:** The data presented in this study are available on request from the corresponding authors.

**Acknowledgments:** Postdoctoral Research Program for Ye. Belyayev, Al-Farabi Kazakh National University, Almaty, Kazakhstan.

**Conflicts of Interest:** The authors declare no conflict of interest.

**Nomenclature**

| | |
|---|---|
| *COP* | Coefficient of Performance |
| $\dot{E}$ | Exergy rate, W |
| *e* | Specific exergy, J/kg |
| *E* | Energy consumption per annum, kWh |
| *GWP* | Global Warming Potential |
| *h* | Specific enthalpy, kJ/kg |
| *L* | Leakage rate per year, %/year |
| $\dot{m}$ | Mass flow rate, kg/s |
| *m* | Refrigerant charge, kg |
| *n* | System operating life, year |
| *P* | Pressure, Pa |
| $\dot{Q}$ | Heat capacity, W |
| *RE* | Relative Exergy |
| *T* | Temperature, °C |
| *TEWI* | Total Equivalent Warming Impact, kg eq. $CO_2$ |
| $\dot{V}$ | Volume flow rate, m³/h |
| *VRC* | Volumetric Refrigeration Capacity, kJ/m³ |
| $\dot{W}$ | Energy consumption, W |

*Greek symbols*

| | |
|---|---|
| *α* | Recovery/recycling factor, % |
| *β* | Indirect $CO_2$ emission factor, kg/kWh |
| *η* | Efficiency |
| *ρ* | Density, kg/m³ |

*Subscripts*

| | |
|---|---|
| EAP | Evaporator Approach |
| elec | Electrical |
| evap | Evaporator, Evaporation |
| CAP | Condenser Approach |
| carnot | Reverse Carnot cycle |
| comp | Compressor |
| cond | Condenser, Condensation |
| cycle | Vapor Compression cycle |
| isen | Isentropic |
| m | Motor |
| pump | Hydraulic loop pump |
| ref | Refrigerant |
| source | Heat source |
| sink | Heat sink |
| SC | Subcooling |
| SH | Superheating |
| system | Overall system |
| t | Transmission |
| TEV | Thermostatic expansion valve |
| total | Total |
| 2nd | The 2nd Law of Thermodynamics |

*Superscripts*

| | |
|---|---|
| alt | Alternative |
| dest | Destruction |
| ex | Exergetic |
| in | Inlet |
| out | Outlet |

*Abbreviation*

| | |
|---|---|
| ATES | Aquifer Thermal Energy Storage |
| BHE | Borehole Heat Exchanger |
| BTES | Borehole Thermal Energy Storage |
| EES | Engineering Equation Solver |
| GEP | Geothermal Energy Pile |
| GHE | Ground Heat Exchanger |
| GSHP | Ground Source Heat Pump |
| HDPE | High-Density Polyethylene |
| HFC | Hydrofluorocarbons |
| HFO | Hydrofluoroolefin |
| HTF | Heat Transfer Fluid |
| TFA | Trifluoroacetate |
| TRT | Thermal Response Test |

**References**

1. Cunha, R.P.; Bourne-Webb, P.J. A critical review on the current knowledge of geothermal energy piles to sustainable climatize buildings. *Renew. Sustain. Energy Rev.* **2022**, *158*, 112072. https://doi.org/10.1016/j.rser.2022.112072.
2. Mohanraj, M.; Belyayev, Ye.; Jayaraj, S.; Kaltayev, A. Research and developments on solar assisted compression heat pump systems-A comprehensive review (Part A: Modeling and modifications). *Renew. Sustain. Energy Rev.* **2018**, *83*, 90–123. https://doi.org/10.1016/j.rser.2017.08.022.
3. Mohanraj, M.; Belyayev, Ye.; Jayaraj, S.; Kaltayev, A. Research and developments on solar assisted compression heat pump systems-A comprehensive review (Part B: Applications). *Renew. Sustain. Energy Rev.* **2018**, *83*, 124–155. https://doi.org/10.1016/j.rser.2017.08.086.
4. Yang, L.W.; Xu, R.J.; Xua, N.; Xia, Y.; Zhou, W.B.; Yang, T.; Belyayev, Ye.; Wang, H.S. Review of the advances in solar-assisted air source heat pumps for the domestic sector. *Energy Convers. Manag.* **2021**, *247*, 114710. https://doi.org/10.1016/j.enconman.2021.114710.
5. Lee, K.C. Classification of geothermal resources–An engineering approach. In Proceedings of the Twenty-First Workshop on Geothermal Reservoir Engineering. Stanford University, Stanford, California, 22–24 January 1996. Available online: https://digital.library.unt.edu/ark:/67531/metadc884944/ (accessed on 3 November 2022).
6. Haehnlein, S.; Beyer, P.; Blum, P. International legal status of the use of shallow geothermal energy. *Renew. Sustain. Energy Rev.* **2010**, *14*, 2611–2625. https://doi.org/10.1016/j.rser.2010.07.069.
7. Saktashova, G.; Aliuly, A.; Belyayev, Ye.; Mohanraj, M.; Singh, R.M. Numerical heat transfer simulation of solar-geothermal hybrid sorce heat pump in Kazakhstan climates. *Bulg. Chem. Commun.* **2018**, *50*, 7–13.
8. Georgiev, A.; Popov, R.; Toshkov, E. Investigation of a hybrid system with ground source heat pump and solar collectors: Charging of thermal storages and space heating. *Renew. Energy* **2020**, *147*, 2774–2790. https://doi.org/10.1016/j.renene.2018.12.087.
9. Vassileva, N.; Georgiev, A.; Popov, R. Simulation study of hybrid ground-source heat pump system with solar collectors. *Bulg. Chem. Commun.* **2016**, *48*, 71–76.
10. Renaldi, R.; Friedrich, D. Techno-economic analysis of a solar district heating system with seasonal thermal storage in the UK. *Appl. Energy* **2019**, *236*, 388–400. https://doi.org/10.1016/j.apenergy.2018.11.030.
11. Welsch, B.; Ruhaak, W.; Schulte, D.O.; Bar, K.; Saas, I. Characteristics of medium deep borehole thermal energy storage. *Int. J. Energy Res.* **2016**, *40*, 855–1868. https://doi.org/10.1002/er.3570.
12. Philippe, M.; Bernier, M.; Marchio, D.; Lopez, S. A semi-analytical model for serpentine horizontal ground heat exchangers. *HVAC&R Res.* **2011**, *17*, 1044–1058.
13. Li, C.; Mao, J.; Peng, X.; Mao, W.; Xing, Z.; Wang, B. Influence of ground surface boundary conditions on horizontal ground source heat pump systems. *Appl. Therm. Eng.* **2019**, *152*, 160–168. https://doi.org/10.1016/j.applthermaleng.2019.02.080.
14. Liu, X.; He, M.; Wang, Y.; Zheng, B.; Li, C.; Zhu, X. Experimental study on heat transfer attenuation due to thermal deformation of horizontal GHEs. *Geothermics* **2021**, *97*, 102241. https://doi.org/10.1016/j.geothermics.2021.102241.
15. Amanzholov, T.; Akhmetov, B.; Georgiev, A.; Kaltayev, A.; Popov, R.; Dzhonova-Atanasova, D.; Tungatarova, M. Numerical modelling as a supplementary tool for Thermal Response Test. *Bulg. Chem. Commun.* **2016**, *48*, 109–114.
16. Akhmetov, B.; Khor, J.O.; Amanzholov, T.; Kaltayev, A.; Romagnoli, A.; Ding, Y. Chapter 12-Modelling at Thermal Energy Storage Device Scale. In *Thermal Energy Storage: Materials, Devices, Systems and Applications*; Ding, Y., Ed.; Royal Society of Chemistry: London, UK, 2021; pp. 370–434. https://doi.org/10.1039/9781788019842-00370.
17. Lee, M.; Lee, D.; Park, M.H.; Kang, Y.T.; Kim, Y. Performance improvement of solar-assisted ground source heat pumps with parallelly connected heat sources in heating-dominated areas. *Energy* **2022**, *240*, 122807. https://doi.org/10.1016/j.energy.2021.122807.
18. Celgia, F.; Marrasso, E.; Roselli, C.; Sasso, M.; Tzscheutschler, P. Exergetic and exergoeconomic analysis of an experimental ground source heat pump system coupled with a thermal storage based on Hardware in Loop. *Appl. Therm. Eng.* **2022**, *212*, 118559. https://doi.org/10.1016/j.applthermaleng.2022.118559.

19. Sang, J.; Liu, X.; Liang, C.; Feng, G.; Li, Z.; Wu, X.; Song, M. Differences between design expectations and actual operation of ground source heat pumps for green buildings in cold region of northern China. *Energy* **2022**, *252*, 124077. https://doi.org/10.1016/j.energy.2022.124077.
20. Hu, R.; Li, X.; Liang, J.; Wang, H.; Liu, G. Field study on cooling performance of a heat recovery ground source heat pump system coupled with thermally activated building systems (TABSs). *Energy Convers. Manag.* **2022**, *262*, 115678. https://doi.org/10.1016/j.enconman.2022.115678.
21. Zhang, M.; Liu, X.; Biswas, K.; Warner, J. A three-dimensional numerical investigation of a novel shallow bore ground heat exchanger integrated with phase change material. *Appl. Therm. Eng.* **2019**, *162*, 114297. https://doi.org/10.1016/j.applthermaleng.2019.114297.
22. Warner, J.; Liu, X.; Shi, L.; Qu, M.; Zhang, M. A novel shallow bore ground heat exchanger for ground source heat pump applications—Model development and validation. *Appl. Therm. Eng.* **2020**, *164*, 114460. https://doi.org/10.1016/j.applthermaleng.2019.114460.
23. Tang, F.; Nowamooz, H. Outlet temperatures of a slinky-type Horizontal Ground Heat Exchanger with the atmosphere-soil interaction. *Renew. Energy* **2020**, *146*, 705–718. https://doi.org/10.1016/j.renene.2019.07.029.
24. Li, W.; Li, X.; Peng, Y.; Wang, Y.; Tu, J. Experimental and numerical studies on the thermal performance of ground heat exchangers in a layered subsurface with groundwater. *Renew. Energy* **2020**, *147*, 620–629. https://doi.org/10.1016/j.renene.2019.09.008.
25. Yang, W.; Xu, R.; Wang, F.; Chen, S. Experimental and numerical investigations on the thermal performance of a horizontal spiral-coil ground heat exchanger. *Renew. Energy* **2020**, *147*, 979–995. https://doi.org/10.1016/j.renene.2019.09.030.
26. Nordell, B.; Snijders, A.L.; Stiles, L. The use of aquifers as thermal energy storage (TES) systems. In *Advances in Thermal Energy Storage Systems, Methods and Applications*; Cabeza, F.C., Ed.; Volume in Woodhead Publishing Series in Energy; Woodhead Publishing: Cambridge, UK, 2015; pp. 87–115.
27. Luo, J.; Pei, K.; Li, P. Analysis of the thermal performance reduction of a groundwater source heat pump (GWHP) system. *Eng. Fail. Anal.* **2022**, *132*, 105922. https://doi.org/10.1016/j.engfailanal.2021.105922.
28. Wang, Z.; Wang, L.; Ma, A.; Liang, K.; Song, Z.; Feng, L. Performance evaluation of ground-water source heat pump system with a fresh air pre-conditioner using ground water. *Energy Convers. Manag.* **2019**, *188*, 250–261. https://doi.org/10.1016/j.enconman.2019.03.061.
29. Boon, D.P.; Farr, G.J.; Abesser, C.; Patton, A.M.; James, D.R.; Schofield, D.I.; Tucker, D.G. Groundwater heat pump feasibility in shallow urban aquifers: Experience from Cardiff, UK. *Sci. Total Environ.* **2019**, *697*, 133847. https://doi.org/10.1016/j.scitotenv.2019.133847.
30. Sanner, B.; Hellstrom, G.; Spitler, J.; Gehlin, S. Thermal response test–current status and world-wide application. In Proceedings of the World Geothermal Congress 2005, Antalya, Turkey, 24–29 April 2005.
31. Feidt, M. *Finite Physical Dimensions Optimal Thermodynamics 1*; Elsevier ISTE Press Ltd.: Amsterdam, The Netherlands, 2017; pp. 245. https://doi.org/10.1016/C2016-0-01231-7.
32. Feidt, M. *Finite Physical Dimensions Optimal Thermodynamics 2*; Elsevier ISTE Press Ltd.: Amsterdam, The Netherlands, 2018; pp. 194. https://doi.org/10.1016/C2017-0-00581-5.
33. Yang, L.W.; Xua, N.; Pu, J.H.; Xia, Y.; Zhou, W.B.; Xu, R.J.; Yang, T.; Belyayev,Ye.; Wang, H.S. Analysis of operation performance of three indirect expansion solar assisted air source heat pumps for domestic heating. *Energy Convers. Manag.* **2022**, *252*, 115061. https://doi.org/10.1016/j.enconman.2021.115061.
34. Yang, L.W.; Li, Y.; Yang, T.; Wang, H.S. Low temperature heating operation performance of a domestic heating system based on indirect expansion solar assisted air source heat pump. *Sol. Energy* **2022**, *244*, 134–154. https://doi.org/10.1016/j.solener.2022.08.037.
35. Colombo, L.P.M.; Lucchini, A.; Molinaroli, L. Experimental analysis of the use of R1234yf and R1234ze(E) as drop-in alternatives of R134a in a water-to-water heat pump. *Int. J. Refrig.* **2020**, *115*, 18–27. https://doi.org/10.1016/j.ijrefrig.2020.03.004.
36. Molinaroli, L.; Lucchini, A.; Colombo, L.P.M. Drop-in analysis of R450A and R513A as low-GWP alternatives to R134a in a water-to-water heat pump. *Int. J. Refrig.* **2022**, *135*, 139–147. https://doi.org/10.1016/j.ijrefrig.2021.12.007.
37. Velasco, F.J.S.; Haddouche, M.R.; Illan-Gomez, F.; Garcia-Cascales, J.R. Experimental characterization of the coupling and heating performance of a $CO_2$ water-to-water heat pump and a water storage tank for domestic hot water production system. *Energy Build.* **2022**, *265*, 112085. https://doi.org/10.1016/j.enbuild.2022.112085.
38. *EES–Engineering Equation Solver*, Professional Version 10.570; F-Chart Software: Madison, WI, USA, 2019.
39. Koshkin, N.N.; Sakun, I.A.; Bambushek, Y.M.; Buharin, N.N.; Gerasimov, Y.D.; Ilyin, A.Y.; Pekarev, V.I.; Stukalenko, A.K.; Timofyevskyi, L.S. *Refrigeration Machines and Installations, Textbook for Technical Colleges*; Mechanical Engineering Publisher (Book in Russian): Saint Petersburg, Russian, 1985; p. 510.
40. Islam, M.A.; Mitra, S.; Thu, K.; Saha, B.B. Study on thermodynamic and environmental effects of vapor compression refrigeration system employing first to next-generation popular refrigerants. *Int. J. Refrig.* **2021**, *131*, 568–580. https://doi.org/10.1016/j.ijrefrig.2021.08.014.
41. Grof, T. *Greening of industry under the Montreal Protocol*; United Nations Industrial Development Organization (UNIDO): Vienna, Austria, 2009; pp. 30.
42. Holland, R.; Khan, M.A.H.; Driscoll, I.; Chantyal-Pun, R.; Derwent, R.G.; Taatjes, C.A.; Orr-Ewing, A.J.; Carl, J.; Percival, C.J.; Shallcross, D.E. Investigation of the production of Trifluoroacetic Acid from two Halocarbons, HFC-134a and HFO-1234yf and its fates using a global three-dimensional chemical transport model. *ACS Earth Space Chem.* **2021,** *5*, 849–857. https://doi.org/10.1021/acsearthspacechem.0c00355.